\documentclass[12pt,a4paper]{article}
\usepackage[dvips]{epsfig}
\usepackage{amsmath}
\usepackage{amssymb}
\usepackage{amsfonts}
\usepackage{amsthm}
\usepackage{amsbsy}
\usepackage{amsgen}
\usepackage{amscd}
\usepackage{amsopn}
\usepackage{amstext}
\usepackage{amsxtra}

\newtheorem{theorem}{Theorem}[section]
\newtheorem{proposition}[theorem]{Proposition}
\newtheorem{lemma}[theorem]{Lemma}
\newtheorem{corollary}[theorem]{Corollary}

\theoremstyle{definition}

\newcommand {\Z} {\mathbb{Z}}
\newcommand {\ZZ} {\mathcal{Z}}

\newcommand {\R} {\mathbb{R}}

\newcommand {\Q} {\mathbb{Q}}

\newcommand {\A} {\mathcal {A}}

\newcommand {\N} {\mathbb {N}}
\newcommand {\M} {\mathcal {M}}
\newcommand {\PP} {\mathbb P}
\newcommand {\OO} {\mathcal {O}}

\begin{document}

\title{Statistics of lattice points in thin annuli for generic lattices}

\author{Igor Wigman \\ School of Mathematical Sciences
\\ Tel Aviv University Tel Aviv 69978, Israel}

\maketitle

\begin{abstract}
We study the statistical properties of the counting function of
lattice points inside thin annuli. By a conjecture of Bleher and
Lebowitz, if the width shrinks to zero, but the area converges to
infinity, the distribution converges to the Gaussian distribution.
If the width shrinks slowly to zero, the conjecture was proven by
Hughes and Rudnick for the standard lattice, and in our previous
paper for generic rectangular lattices. We prove this conjecture
for arbitrary lattices satisfying some generic Diophantine
properties, again assuming the width of the annuli shrinks slowly
to zero.

One of the obstacles of applying the technique of Hughes-Rudnick
on this problem is the existence of so-called close pairs of
lattice points. In order to overcome this difficulty, we bound the
rate of occurence of this phenomenon by extending some of the work
of Eskin-Margulis-Mozes on the quantitative Openheim conjecture.

\end{abstract}

\section{Introduction}

We consider a variant of the lattice points counting problem. Let
$\Lambda \subset \R^2$ be a planar lattice, with $\det\Lambda$ the
area of its fundamental cell.  Let
\begin{equation*}
N_{\Lambda} (t) = \{ x\in\Lambda:\: |x| \le t\},
\end{equation*}
denote its counting function, that is, we are counting
$\Lambda$-points inside a disc of radius $t$.

As well known, as $t\rightarrow\infty$, $N_{\Lambda} (t) \sim
\frac{\pi}{\det{\Lambda}}t^2$. Denoting the remainder or the error
term
\begin{equation*}
\Delta_{\Lambda} (t) = N_{\Lambda} (t) - \frac{\pi}{\det{\Lambda}}
t^2,
\end{equation*}
it is a conjecture of Hardy that
\begin{equation*}
|\Delta_{\Lambda} (t)| \ll_{\epsilon} t^{1/2+\epsilon}.
\end{equation*}

Another problem one could study is the {\em statistical} behavior
of the value distribution of $\Delta_{\Lambda}$ normalized by
$\sqrt{t}$, namely of
\begin{equation*}
F_{\Lambda} (t) := \frac{\Delta_{\Lambda} (t)}{\sqrt{t}}.
\end{equation*}

Heath-Brown ~\cite{HB} shows that for the standard lattice
$\Lambda = \Z^2$, the value distribution of $F_{\Lambda}$, weakly
converges to a non-Gaussian distribution with density $p(x)$.
Bleher ~\cite{BL3} established an analogue of this theorem for a
more general setting, where in particular it implies a
non-Gaussian limiting distribution of $F_{\Lambda}$, for any
lattice $\Lambda\subset\Z^2$.

However, the object of our interest is slightly different. Rather
than counting lattice points in the circle of varying radius $t$,
we will do the same for {\em annuli}. More precisely, we define
\begin{equation*}
N_{\Lambda} (t, \, \rho) := N_{\Lambda} (t+\rho) - N_{\Lambda}
(t),
\end{equation*}
that is, the number of $\Lambda$-points inside the annulus of
inner radius $t$ and width $\rho$. The "expected" value is the
area $\frac{\pi}{\det{\Lambda}} (2t\rho+\rho^2)$, and the
corresponding normalized remainder term is
\begin{equation*}
S_{\Lambda} (t,\,\rho ) := \frac{N_{\Lambda} (t+\rho) -
N_{\Lambda} (t) - \frac{\pi}{\det{\Lambda}} (2t\rho+\rho^2)
}{\sqrt{t}}.
\end{equation*}

The statistics of $S_{\Lambda} (t,\,\rho )$ vary depending to the
size of $\rho(t)$. Of our particular interest is the {\em
intermediate} or {\em macroscopic regime}. Here $\rho\rightarrow
0$, but $\rho t \rightarrow\infty$. A particular case of the
conjecture of Bleher and Lebowitz \cite{BL4} states that
$S_{\Lambda} (t,\,\rho)$ has a Gaussian distribution. In 2004
Hughes and Rudnick ~\cite{HR} established the Gaussian
distribution for the unit circle, under an additional assumption
that $\rho(t) \gg t^{-\epsilon}$ for every $\epsilon > 0$.

By a rotation and dilation (which does not essentially effect the
counting function), we may assume, with no loss of generality,
that $\Lambda$ admits a basis one of whose elements is the vector
$(1,0)$, that is $\Lambda = \big\langle 1, \,\alpha+i\beta
\big\rangle$ (we make the natural identification of $i$ with
$(0,\, 1)$). In a previous paper ~\cite{W} we already dealt with
the problem of investigating the statistical properties of the
error term for rectangular lattice $\Lambda=\big\langle 1,\,
i\beta\big\rangle$. We established the limiting Gaussian
distribution for the "generic" case in this 1-parameter family.

Some of the work done in ~\cite{W} extends quite naturally for the
2-parameter family of planar lattices $\big\langle 1, \,
\alpha+i\beta \big\rangle$. That is, in the current work we will
require algebraic independence of $\alpha$ and $\beta$ as well as
the "strong Diophantinity" of the \underline{pair} $(\alpha,\,
\beta)$ (to be defined), rather than transcendence and strong
Diopantinity of the aspect ratio of the ellipse, as in ~\cite{W}.

We say that a real number $\xi$ is {\em strongly Diophantine}, if
for every {\em fixed} natural $n$, there exists $K_1>0$, such that
for integers $a_j$ with $\sum\limits_{j=0}^{n} a_j \alpha^j \ne
0$,
\begin{equation*}
\bigg| \sum\limits_{j=0}^{n} a_j \xi^j \bigg| \gg_{n}
\frac{1}{\bigg( \max\limits_{0\le j \le n} |a_j| \bigg) ^ {K_1}}.
\end{equation*}
It was shown by Mahler ~\cite{MAH}, that this property holds for a
"generic" real number.
We say that a pair of numbers $(\alpha, \, \beta)$ is {\em
strongly Diophantine}, if for every {\em fixed} natural $n$, there
exists a number $K_1 > 0$, such that for every integral polynomial
$p(x,\, y) = \sum\limits_{i+j \le n} {} a_{i,\, j} x^i y^j$ of
degree $\le n$, we have
\begin{equation*}
| p(\alpha,\, \beta)| \gg_{n} \frac{1}{\max\limits_{i+j \le n}
|a_{i,\, j}| ^ {K_1}},
\end{equation*}
whenever $p(\alpha,\beta)\neq 0$. 
This holds for almost all real pairs $(\alpha, \, \beta)$, see
section~\ref{subsec:high mom}.

\begin{theorem}
\label{thm:norm dist} Let $\Lambda = \big\langle 1,
\,\alpha+i\beta \big\rangle$ where $(\alpha, \, \beta)$ is
algebraically independent and strongly Diophantine pair of real
numbers. Assume that $\rho = \rho (T) \rightarrow 0$, but for
every $\delta > 0$, $\rho \gg T^{-\delta}$. Then for every
interval $\A$,
\begin{equation}
\lim_{T\rightarrow \infty } \frac{1}{T} meas  \bigg\{ t\in [T,\,
2T ] :\: \frac{S_{\Lambda} (t,\, \rho ) } {\sigma } \in \A \bigg\}
= \frac{1}{\sqrt{2\pi }} \int\limits_{\A} e ^ {-\frac{x^2}{2}} dx,
\end{equation}
where the variance is given by
\begin{equation}
\label{eq:sigma comp rho A} \sigma^2 := \frac{4 \pi}{\beta} \cdot
\rho.
\end{equation}
\end{theorem}

\paragraph{Remark:} Note that the variance $\sigma^2$ is
$\alpha$-independent, since the determinant $\det (\Lambda) =
\beta$.

One of the features of a rectangular lattice is that it is quite
easy to show that the number of so-called close pairs of lattice
points or pairs of points lying within a narrow annulus is bounded
by essentially its average (see lemma 5.2 of ~\cite{W}). This
particular feature of the rectangular lattices was exploited while
reducing the computation of the moments to the ones of a smooth
counting function (we call it "unsmoothing"). In order to prove an
analogous bound for a general lattice, we extend a result from
Eskin, Margulis and Mozes ~\cite{EMM} for our needs to obtain
proposition \ref{prop:S_n bnd blh}. We believe that this
proposition is of independent interest.

\section{The distribution of $\tilde{S}_{\Lambda , \, M , \, L}$ }
\label{sec:smth dist}

In this section, we are interested in the distribution of the
smooth version of $S_{\Lambda} (t,\, \rho )$, denoted
$\tilde{S}_{\Lambda , \, M, \, L} (t)$, where $L :=
\frac{1}{\rho}$ and $M$ is the smoothing parameter. Just as in
~\cite{W} and ~\cite{HR},
\begin{equation}
\label{eq:ann cnt smth} \tilde{S}_{\Lambda , \, M, \, L} (t) =
\frac { \tilde{N}_{\Lambda ,\, M} (t + \frac{1}{L} ) -
                                \tilde{N}_{\Lambda ,\, M}
                                (t) - \frac{\pi}{d} (\frac{2 t}{L} +
\frac{1}{L^{2}})}{\sqrt {t} },
\end{equation}
where $\tilde{N}_{\Lambda ,\, M}$ is the smooth version of
$N_{\Lambda}$, computed by means of convolution of the
characteristic function of the unit ball with $\psi$, a smooth
function with a compact support (see ~\cite{HR} or ~\cite{W} for
details). We assume that for every $\delta
> 0$, $L = L(T) = O(T^{\delta } )$, which corresponds to the
assumption of theorem \ref{thm:norm dist} regarding $\rho :=
\frac{1}{L}$.

Rather than drawing $t$ at random from $[T,\, 2T]$ with a uniform
distribution, we prefer to work with smooth densities: introduce
$\omega \ge 0$, a smooth function of total mass unity, such that
both $\omega$ and $\hat{\omega}$ are rapidly decaying, namely
\begin{equation*}
|\omega(t) | \ll \frac{1}{(1+|t|)^A} , \; \; | \hat{\omega}(t) |
\ll \frac{1}{(1+|t|)^A},
\end{equation*}
for every $A>0$. Define the averaging operator
\begin{equation*}
\langle f \rangle _T = \frac{1}{T} \int\limits_{-\infty}^{\infty}
f(t) \omega (\frac{t}{T} ) dt ,
\end{equation*}
and let $\PP_{\omega , \, T }$ be the associated probability
measure:
\begin{equation*}
\PP_{\omega , \, T } (f \in \A ) = \frac{1}{T}
\int\limits_{-\infty}^{\infty} 1_{\A} (f(t)) \omega (\frac{t}{T} )
dt.
\end{equation*}
\paragraph{Remark:} In what follows, we will suppress the explicit
dependency on $T$, whenever convenient.

\begin{theorem}
\label{thm:norm dist smth} Suppose that $M(T)$ and $L(T)$ are
increasing to infinity with $T$, such that $M=O(T^{\delta})$ for
all $\delta > 0$, and $L/\sqrt{M} \rightarrow 0$. Then if
$(\alpha,\, \beta)$ is an {\em algebraically independent strongly
Diophantine} pair, we have for $\Lambda = \big\langle 1, \, \alpha
+ i \beta \big\rangle$,
\begin{equation*}
\lim_{T\rightarrow \infty } \PP_{\omega,\, T}  \bigg\{
\frac{\tilde{S}_{\Lambda, \, M, \, L} } {\sigma } \in \A \bigg\} =
\frac{1}{\sqrt{2\pi }} \int\limits_{\A} e ^ {-\frac{x^2}{2}} dx,
\end{equation*}
for any interval $\A$, where
\begin{equation}
\label{eq:sigma comp rho beta} \sigma^2 := \frac{4 \pi }{\beta L}.
\end{equation}
\end{theorem}

\paragraph{Definition:} A tuple of real numbers
$(\alpha_1,\, \ldots,\, \alpha_n ) \in \R^n$ is called {\em
Diophantine}, if there exists a number $K > 0$, such that for
every integer tuple $\{ a_i \} _{i=0}^{n}$,
\begin{equation}
\label{eq:dioph num} \bigg|  a_0 + \sum\limits_{i=1}^{n} a_i
\alpha_i \bigg| \gg \frac{1}{q^{K}},
\end{equation}
where $q = \max\limits_{0\le i \le n} {|a_i|}$. Khintchine proved
that {\em almost all} tuples in $\R^n$ are Diophantine (see, e.g.
~\cite{S}, pages 60-63).

Denote the dual lattice
\begin{equation*}
\Lambda ^{*} = \big\langle 1,\, -\frac{\alpha}{\beta} + i
\frac{1}{\beta} \big\rangle.
\end{equation*}
We assume for the rest of current section that the set of squared
norms of $\Lambda ^{*}$ satisfy the Diophantine property, which
means that $(\alpha ^2, \, \alpha\beta,\, \beta^2 )$ is a
Diophantine triple of numbers. We may assume the Diophantinity of
$(\alpha ^2, \, \alpha\beta,\, \beta^2 )$, since theorem
\ref{thm:norm dist} (and theorem \ref{thm:norm dist smth}) assume
$(\alpha,\, \beta)$ is {\em strongly Diophantine}, which is
obviously a stronger assumption.

We use the following approximation to $\tilde{N}_{\Lambda, M} (t)$
(see e.g ~\cite{W}, lemma 4.1):
\begin{lemma}
As $t\rightarrow \infty$,
\begin{equation}
\label{eq:app smth cnt} \tilde{N}_{\Lambda, M} (t) = \frac{\pi
t^2}{\beta} - \frac{\sqrt {t}}{\beta \pi} \sum\limits_{\vec{k} \in
\Lambda ^ {*} \setminus \{ 0 \} } \frac{\cos \big( 2\pi t |\vec{k}
| + \frac{\pi}{4} \big) } {|\vec{k}|^{\frac{3}{2}}} \cdot
        \hat{\psi} \bigg( \frac{|\vec{k}|}{\sqrt{M}} \bigg) +
O\bigg(\frac{1}{\sqrt {t}} \bigg),
\end{equation}
where, again, $\Lambda^*$ is the dual lattice.
\end{lemma}

By the definition of $\tilde{S}_{\Lambda, \, M ,\, L }$ in
\eqref{eq:ann cnt smth} and appropriately manipulating the sum in
\eqref{eq:app smth cnt} we obtain the following
\begin{corollary}
\begin{equation}
\begin{split}
\label{eq:approx S smth} \tilde{S}_{\Lambda, \, M ,\, L } (t) &=
\frac{2}{\beta \pi} \sum\limits_{\vec{k} \in \Lambda ^ {*}
\setminus \{ 0 \} } \frac{\sin \bigg( \frac{\pi |\vec{k}|}{L}
\bigg) } {|\vec{k}| ^{\frac{3}{2}}}
 \sin \bigg( 2 \pi \big( t+\frac{1}{2 L} \big) |\vec{k} | + \frac{\pi}{4} \bigg)
\hat{\psi} \bigg( \frac{|\vec{k}|}{\sqrt {M}} \bigg)  \\
&+ O \bigg(\frac{1}{\sqrt{t}} \bigg)\;.
\end{split}
\end{equation}
\end{corollary}

One should note that $\hat{\psi}$ being compactly supported means
that the sum essentially truncates at $|\vec{k}| \approx
\sqrt{M}$.

Unlike the standard lattice, clearly there are no nontrivial
multiplicities in $\Lambda$, that is
\begin{lemma}
\label{lem:no mult} Let $\vec{a_j} = m_j + n_j (\alpha+i\beta) \in
\Lambda$, $j = 1,\, 2$, with an irrational $\alpha$ such that
$\gamma\notin\Q(\alpha)$. Then if $|\vec{a_1}| = |\vec{a_2}|$,
either $n_1 = n_2$ and $m_1=m_2$ or $n_1 = -n_2$ and $n_2 = -m_2$.
\end{lemma}

\begin{proof}[Proof of theorem \ref{thm:norm dist smth}]
We will show that the moments of $\tilde{S}_{\Lambda, \, M ,\, L
}$ corresponding to the smooth probability space converge to the
moments of the normal distribution with zero mean and variance
which is given by theorem \ref{thm:norm dist smth}. This allows us
to deduce that the distribution of $\tilde{S}_{\Lambda, \, M ,\, L
}$ converges to the normal distribution as $T\rightarrow\infty$,
precisely in the sense of theorem \ref{thm:norm dist smth}.

First, we show that the mean is $O(\frac{1}{\sqrt{T}})$. Since
$\omega$ is real,
\begin{equation*}
\Bigg| \Bigg\langle \sin \bigg( 2 \pi \big( t+\frac{1}{2 L} \big)
|\vec{k} | + \frac{\pi}{4} \bigg) \Bigg\rangle \Bigg| = \bigg| \Im
m \bigg\{ \hat{\omega} \big (-T |\vec{k}| \big) e ^ { i \pi
(\frac{|\vec{k}|}{L} + \frac{1}{4} } \bigg\} \bigg| \ll
\frac{1}{T^A |\vec{k}| ^ {A}}
\end{equation*}
for any $A>0$, where we have used the rapid decay of
$\hat{\omega}$. Thus
\begin{equation*}
\bigg| \bigg\langle \tilde{S}_{\Lambda, \, M ,\, L } \bigg\rangle
\bigg| \ll \sum\limits_{\vec{k} \in \Lambda ^ {*} \setminus \{ 0
\} } \frac{1}{T^A |\vec{k}|^ {A+3/2}} + O \bigg(
\frac{1}{\sqrt{T}} \bigg) \ll  O \bigg( \frac{1}{\sqrt{T}} \bigg),
\end{equation*}
due to the convergence of $\sum\limits_{\vec{k} \in \Lambda ^ {*}
\setminus \{ 0 \} } \frac{1}{|\vec{k}|^ {A+3/2}}$, for $A >
\frac{1}{2}$

Now define
\begin{equation}
\label{eq:def M} \M_{\Lambda, \, m} := \Bigg\langle \bigg(
\frac{2}{\beta \pi} \sum\limits_{\vec{k} \in \Lambda ^ {*}
\setminus \{ 0 \} } \frac{\sin \bigg( \frac{\pi |\vec{k}|}{L}
\bigg) } {|\vec{k}| ^{\frac{3}{2}}}
 \sin \bigg( 2 \pi \big( t+\frac{1}{2 L} \big) |\vec{k} | + \frac{\pi}{4} \bigg)
\hat{\psi} \big( \frac{|\vec{k}|}{\sqrt {M}} \big) \bigg) ^{m}
\Bigg\rangle
\end{equation}

Then from \eqref{eq:approx S smth}, the binomial formula and the
Cauchy-Schwartz inequality,
\begin{equation*}
\bigg\langle \big( \tilde{S}_{\Lambda, \, M ,\, L }  \big) ^m
\bigg\rangle = \M_{\Lambda, \, m} + O \bigg( \sum
\limits_{j=1}^{m} \binom{m}{j} \frac{\sqrt {\M_{2m-2j}} }{T^{j/2}}
 \bigg)
\end{equation*}

Proposition \ref{prop:var comp} together with proposition
\ref{prop:high mom comp} allow us to deduce the result of theorem
\ref{thm:norm dist smth} for an algebraically independent strongly
Diophantine $(\xi,\, \eta) := (-\frac{\alpha}{\beta},\,
\frac{1}{\beta})$. Clearly, $(\alpha,\, \beta)$ being
algebraically independent and strongly Diophantine is sufficient.
\end{proof}

\subsection{The variance}
\label{ssec:var comp} The computation of the variance is done in
two steps. First, we reduce the main contribution to the {\em
diagonal} terms, using the assumption on the pair $(\alpha,\,
\beta)$ (i.e. $(\alpha^2,\, \alpha\beta,\, \beta^2)$ is {\em
Diophantine}). Then we compute the contribution of the {\em
diagonal} terms. We sketch these steps, since they are very close
to the corresponding one ~\cite{W}.

Suppose that the triple $(\alpha ^2, \, \alpha\beta,\, \beta^2 )$
satisfies \eqref{eq:dioph num}.
\begin{proposition}
\label{prop:var comp} If $M = O \big( T^{1/(K + 1/2 +\delta)}
\big)$ for fixed $\delta > 0$, then the variance of
$\tilde{S}_{\Lambda, \, M ,\, L } $ is asymptotic to
\begin{equation*}
\sigma^2 :=
 \frac{4}{\beta^2 \pi^2}
\sum\limits_{\vec{k} \in \Lambda ^ {*} \setminus \{ 0 \} }
\frac{\sin ^2 \bigg( \frac{\pi |\vec{k}|}{L} \bigg) } {|\vec{k}|
^{3}} \hat{\psi} ^ 2 \bigg( \frac{|\vec{k}|}{\sqrt {M}} \bigg)
\end{equation*}
If $L \rightarrow \infty$, but $L/ \sqrt {M} \rightarrow 0$, then
\begin{equation}
\label{eq:sigma comp} \sigma^2 \sim \frac{4 \pi }{\beta L}
\end{equation}
\end{proposition}

\begin{proof}
Expanding out \eqref{eq:def M}, we have
\begin{equation}
\label{eq:M exp}
\begin{split}
\M_{\Lambda, \, 2} =  \frac{4}{\beta^2 \pi ^2} &
\sum\limits_{\vec{k},\vec{l} \in \Lambda ^ {*} \setminus \{ 0 \} }
 \frac{\sin \bigg( \frac{\pi |\vec{k}|}{L} \bigg)
\sin \bigg( \frac{\pi |\vec{l}|}{L} \bigg) \hat{\psi} \big(
\frac{|\vec{k}|}{\sqrt {M}} \big) \hat{\psi} \big(
\frac{|\vec{l}|}{\sqrt {M}} \big) }
{|\vec{k}| ^{\frac{3}{2}} |\vec{l}| ^{\frac{3}{2}} } \\
& \times \bigg\langle \sin \bigg( 2 \pi \bigg( t+\frac{1}{2 L}
\bigg) |\vec{k} | + \frac{\pi}{4} \bigg) \sin \bigg( 2 \pi \bigg(
t+\frac{1}{2 L} \bigg) |\vec{l} | + \frac{\pi}{4} \bigg)
\bigg\rangle
\end{split}
\end{equation}
It is easy to check that the average of the second line of the
previous equation is:
\begin{equation}
\label{eq:mean sines}
\begin{split}
\frac{1}{4}\bigg[ &\hat{\omega} \big( T(|\vec{k}| - |\vec{l}|)
\big) e^{i\pi (1/L) (|\vec{l}| -
|\vec{k}|)} + \\
&\hat{\omega} \big( T(|\vec{l}| - |\vec{k}|) \big) e^{i\pi (1/L)
(|\vec{k}| -
|\vec{l}|)} + \\
&\hat{\omega} \big( T(|\vec{k}| + |\vec{l}|) \big) e^{-i\pi (1/2 +
(1/L)
(|\vec{k}| + |\vec{l}|))} - \\
&\hat{\omega} \big( -T(|\vec{k}| + |\vec{l}|) \big) e^{i\pi (1/2 +
(1/L) (|\vec{k}| + |\vec{l}|))} \bigg]
\end{split}
\end{equation}
Recall that the support condition on $\hat{\psi}$ means that
$\vec{k}$ and $\vec{l}$ are both constrained to be of length
$O(\sqrt{M} )$. Thus the off-diagonal contribution (that is for
$|\vec{k}| \ne |\vec{l}|$ ) of the first two lines of
\eqref{eq:mean sines} is
\begin{equation*}
\ll \sum\limits_ {\substack {\vec{k},\vec{l} \in \Lambda ^ {*}
\setminus \{ 0 \} \\ |\vec{k}|,\,|\vec{k'}| \le \sqrt {M}}}
\frac{M^{A(K+1/2)}}{T^A} \ll \frac{M^{A(K+1/2) + 2}}{T^A} \ll
T^{-B},
\end{equation*}
for every $B>0$, by Diophantinity of $(\alpha,\, \alpha\beta, \,
\beta^2)$.

Obviously, the contribution to \eqref{eq:M exp} of the two last
lines of \eqref{eq:mean sines} is negligible both in the diagonal
and off-diagonal cases, justifying the diagonal approximation of
\eqref{eq:M exp} in the first statement of the proposition, and we
omit the rest of the proof.
\end{proof}

\subsection{The higher moments}
\label{subsec:high mom} In order to compute the higher moments we
will prove that the main contribution comes from the so-called
{\em diagonal} terms (to be explained later). In order to bound
the contribution of the {\em off-diagonal} terms, we will use the
following theorem, which is a consequence of the work of Kleinbock
and Margulis ~\cite{KM}. The contribution of the diagonal terms is
computed exactly in the same manner it was done in ~\cite{W}, and
so we will omit it here. 

\begin{theorem}
\label{thm:alm all str dphnt} Let an integer $n$ be given. Then
almost all pairs of real numbers $(\xi,\,\eta)\in \R^2$ satisfy
the following property: there exists a number $K_1 \in \N$ such
that for every integer polynomial of $2$ variables $p(x,\, y) =
\sum\limits_{i+j \le n} {} a_{i,\, j} x^i y^j$ with degree $\le
n$, we have
\begin{equation*}
\big| p(\xi,\, \eta ) \big| \gg h^{-K_1},
\end{equation*}
where $h=\max\limits_{i+j\le n} {|a_{i,\, j}|}$ is the height of
$p$. The constant involved in the $"\gg"$ notation depends only on
$\xi,\, \eta$ and $n$.
\end{theorem}

We will remark that theorem A in ~\cite{KM} is much more general
when the result we are using. As a matter of fact, we have the
inequality
\begin{equation*}
 \big| b_0 + b_1 f_1 (x) + \ldots + b_n f_n (x) \big| \gg_{\epsilon}
\frac{1}{h^{n + \epsilon}}
\end{equation*}
with $b_i\in\Z$ and
\begin{equation*}
h := \max\limits_{0 \le i \le n} |b_i|.
\end{equation*}
The inequality above holds for every $\epsilon > 0$ for a wide
class of functions $f_i : U\rightarrow\R$, for almost all $x\in
U$, where $U\subset \R^{m}$ is an open subset. Here we use this
inequality for the monomials.

\paragraph{Definition:} We call the pairs $(\xi,\, \eta)$ which satisfy for all
natural $n$ the property of theorem \ref{thm:alm all str dphnt},
{\em strongly Diophantine}. Thus theorem \ref{thm:alm all str
dphnt} states that almost all real pairs of numbers are strongly
Diophantine.

\paragraph{Remark:} Simon Kristensen ~\cite{KR} has recently shown,
that the set of all pairs $(\xi,\, \eta)\in\R^2$ which fail to be
strongly Diophantine has Hausdorff dimension $1$.

Obviously, strong Diophantinity of $(\xi,\, \eta)$ implies
Diophantinity of any $n$-tuple of real numbers which consists of
any set of monomials in $\xi$ and $\eta$. Moreover, $(\xi,\,\eta)$
is strongly Diophantine iff $(-\frac{\alpha}{\beta},\,
\frac{1}{\beta})$ is such.


We have the following analogue of lemma 4.7 in ~\cite{W}, which
will eventually allow us to exploit the strong Diophantinity of
$(\alpha,\,\beta)$.
\begin{lemma}
\label{lem:gen bnd eps} If $(\xi,\, \eta)$ is strongly
Diophantine, then it satisfies the following property: for any
fixed natural $m$, there exists $K \in \N$, such that if
\begin{equation*}
z_j = a_j^2+b_j^2\xi^2 + 2a_j b_j\xi +  b_j^2 \eta^2\ll M,
\end{equation*}
and $\epsilon_j = \pm 1$ for $j=1,\ldots , m$, with integral
$a_j,\, b_j$ and if $\sum\limits_{j=1}^{m} \epsilon_j \sqrt{z_j}
\ne 0$, then
\begin{equation}
\label{eq:gen bnd res} \big| \sum\limits_{j=1}^{m} \epsilon_j
\sqrt{z_j} \big| \gg M ^ {-K},
\end{equation}
where the constant involved in the $"\gg"$ notation depends only
on $\eta$ and $m$.
\end{lemma}

The proof is essentially the same as the one of lemma 4.7 from
~\cite{W}, considering the product $Q$ of numbers of the form
$\sum\limits_{j=1}^{m} \delta_j \sqrt{z_j}$ over all possible
signs $\delta_j$. Here we use the Diophantinity of the real tuple
$(\xi,\, \eta)$ rather than of a single real number.

\begin{proposition}
\label{prop:high mom comp} Let $m\in\N$ be given. Suppose that
$\Lambda = \langle 1, \, \alpha + i \beta \rangle$, such that the
pair $(\xi,\, \eta) := (-\frac{\alpha}{\beta},\, \frac{1}{\beta})$
is algebraically independent strongly Diophantine, which satisfy
the property of lemma \ref{lem:gen bnd eps} for the given $m$,
with $K=K_m$. Then if $\M = O \big( T^{\frac{1-\delta}{K_m}}
\big)$ for some $\delta > 0$, and if $L \rightarrow \infty$ such
that $L / \sqrt{M} \rightarrow 0$, the following holds:
\begin{equation*}
\frac{\M_{\Lambda, \, m} }{\sigma^m} =
\begin{cases}
\frac{m!}{2^{m/2} \big( \frac{m}{2} \big) !  } + O\big( \frac{\log
L}{L} \big),
\; &m \text{ is even} \\
O\big( \frac{\log L}{L} \big), \; &m \text{ is odd}
\end{cases}
\end{equation*}
\end{proposition}

\begin{proof}
Expanding out \eqref{eq:def M},  we have
\begin{equation}
\label{eq:exp M_m}
\begin{split}
\M_{\Lambda, \, m} =  \frac{2^m}{\beta^m \pi ^m}
 \sum\limits_{\vec{k_1},\ldots ,\,\vec{k_m} \in \Lambda ^ {*} \setminus \{ 0 \}
} &\prod\limits_{j=1}^{m}
 \frac{\sin \bigg( \frac{\pi |\vec{k_j}|}{L} \bigg)
\hat{\psi} \big( \frac{|\vec{k_j}|}{\sqrt {M}} \big)}
{|\vec{k_j}| ^{\frac{3}{2}}} \\
& \times \bigg\langle \prod\limits_{j=1}^{m} \sin \bigg( 2 \pi
\big( t+\frac{1}{2 L} \big) |\vec{k_1} | + \frac{\pi}{4} \bigg)
\bigg\rangle
\end{split}
\end{equation}

Now,
\begin{equation*}
\begin{split}
\bigg\langle \prod\limits_{j=1}^{m} &\sin \bigg( 2 \pi \big(
t+\frac{1}{2 L} \big) |\vec{k_1} | + \frac{\pi}{4}
\bigg) \bigg\rangle \\
&= \sum \limits _{\epsilon_j = \pm 1} \frac{\prod\limits_{j=1}^{m}
\epsilon_j} {2^m i^m} \hat{\omega} \bigg( -T \sum\limits_{j=1}^{m}
\epsilon_j |\vec{k_j}| \bigg) e^{ \pi i \sum\limits_{j=1} ^{m}
\epsilon_j \big( (1/L) |\vec{k_j}| + 1/4 \big) }
\end{split}
\end{equation*}

We call a term of the summation in \eqref{eq:exp M_m} with
$\sum\limits_{j=1}^{m} \epsilon_j |\vec{k_j}| = 0$ {\em diagonal},
and {\em off-diagonal} otherwise. Due to lemma \ref{lem:gen bnd
eps}, the contribution of the {\em off-diagonal} terms is:

\begin{equation*}
\ll \sum\limits_{\vec{k_1},\ldots ,\,\vec{k_m} \in \Lambda ^ {*}
\setminus \{ 0 \} } \bigg( \frac{T}{M ^ {K_m}} \bigg) ^ {-A} \ll
M^m T^{-A\delta},
\end{equation*}
for every $A>0$, by the rapid decay of $\hat{\omega}$ and our
assumption regarding $M$.

Since $m$ is constant, this allows us to reduce the sum to the
{\em diagonal terms}. In order to be able to sum over all the
diagonal terms we need the following analogue of a well-known
theorem due to Besicovitch ~\cite{BS} about incommensurability of
square roots of integers.

\begin{proposition}
\label{prop:diag impl prnc diag} Suppose that $\xi$ and $\eta$ are
algebraically independent, and
\begin{equation}
\label{eq:norms expr} z_j = a_j^2+2 a_j b_j \xi + b_j^2
(\xi^2+\eta),
\end{equation}
such that $(a_j,\, b_j ) \in \Z_+ ^2$ are all different primitive
vectors, for $1\le j \le m$. Then $\{ \sqrt{z_j} \} _{j=1} ^{m} $
are linearly independent over $\Q$.
\end{proposition}

The last proposition is an immediate consequence of a theorem
proved in the appendix of ~\cite{BL2}.

Computing the contribution of the diagonal terms is done literally
the same way it was done in ~\cite{W} and thus it is omitted here.
In order to be able to sum over the diagonal terms, we use here
proposition \ref{prop:diag impl prnc diag}, rather than
proposition 3.2 in ~\cite{W}.

\end{proof}

\section{Bounding the number of close pairs of lattice points}
\label{sec:bnd cls prs} Roughly speaking, we say that a pair of
lattice points, $n$ and $n'$ is {\em close}, if $\big| |n|-|n'|
\big|$ is {\em small}. We would like to show that this phenomenon
is {\em rare}. This is closely related to the Oppenheim
conjecture, as $|n|^2-|n'|^2$ is a quadratic form on the
coefficients of $n$ and $n'$. In order to establish a quantative
result, we use a technique developed in a paper by Eskin, Margulis
and Mozes ~\cite{EMM}.

\subsection{Statement of the results}
The ultimate goal of this section is to establish the following
\begin{proposition}
\label{prop:S_n bnd blh} Let $\Lambda$ be a lattice and denote
\begin{equation}
\label{eq:A(R,delta) def} A(R,\delta) := \{ (\vec{k},\, \vec{l})
\in \Lambda :\: R \le |\vec{k}|^2 \le 2R,\, |\vec{k}|^2 \le
|\vec{l}|^2 \le |\vec{k}|^2 + \delta \}.
\end{equation}
Then if $\delta > 1$, such that $\delta = o(R)$, we have
\begin{equation*}
\# A(R,\delta) \ll R \delta \cdot \log{R}
\end{equation*}
\end{proposition}

In order to prove this result, we note that evaluating the size of
$A(R,\, \delta)$ is equivalent to counting integer points
$\vec{v}\in\R^4$ with $T \le \| \vec{v} \| \le 2T$ such that
\begin{equation*}
0 \le Q_{1}(v) \le \delta,
\end{equation*}
where $Q_{1}$ is a quadratic form of signature $(2,\, 2)$, given
explicitly by
\begin{equation}
\label{eq:quad form def} Q_{1} (\vec{v}) = (v_1 + v_2 \alpha)^2 +
(v_2 \beta)^2 - (v_3 + v_4 \alpha)^2 - (v_4 \beta)^2.
\end{equation}
For a fixed $\delta>0$ and a large $R$, this situation was
considered extensively by Eskin, Margulis and Mozes ~\cite{EMM}.
We will examine how the constants involved in their result depend
on $\delta$, and find out that there is a linear dependency, which
is what we essentially need. The author wishes to thank Alex Eskin
for his assistance with this matter.

\paragraph{Remark:} For our purposes we need a weaker result:
\begin{equation*}
\# A(R,\delta) \ll_{\epsilon} R \delta \cdot R^{\epsilon},
\end{equation*}
for every $\epsilon > 0$. If $\Lambda$ is a rectangular lattice
(i.e. $\alpha = 0$), then this result follows from properties of
the divisor function (see e.g. ~\cite{BL}, lemma 3.2).

Theorem 2.3 in ~\cite{EMM} considers a more general setting than
proposition \ref{prop:S_n bnd blh}. We state here theorem 2.3 from
~\cite{EMM} (see theorem \ref{thm:mem thm 2.3}). It follows from
theorem 3.3 from ~\cite{EMM}, which will be stated as well (see
theorem \ref{thm:mem thm 3.3}). Then we give an outline of the
proof of theorem 2.3 of ~\cite{EMM}, and inspect the dependency on
$\delta$ of the constants involved.

\subsection{Theorems 2.3 and 3.3 from ~\cite{EMM}} Let $\Delta$
be a lattice in $\R^n$. We say that a subspace $L\subset\R^n$ is
$\Delta$-rational, if $L\cap\Delta$ is a lattice in $L$. We need
the following definitions:
\paragraph{Definitions:}
\begin{equation*}
\alpha_ i (\Delta) := \sup\bigg\{ \frac{1}{d_{\Delta} (L) } \bigg|
\: L \text{ is a } \Delta-\text{rational subspace of dimension } i
\bigg\},
\end{equation*}
where
\begin{equation*}
d_{\Delta} (L) := vol (L / (L \cap \Delta) ).
\end{equation*}
Also
\begin{equation*}
\alpha (\Delta) := \max\limits_{0\le i \le n} \alpha_i (\Delta).
\end{equation*}

Since the space of unimodular lattices is canonically isomorphic
to \\ $SL(n,\, \R) / SL(n,\, \Z)$, the notation $\alpha (g)$ makes
sense for $g\in G:=SL(n,\, \R)$.

For a bounded function $f : \R^n\rightarrow \R$, with $|f|\le M$,
which vanishes outside a ball $B(0,\, R)$, define
$\tilde{f}:SL(n,\, \R) \rightarrow \R$ by the following formula:
\begin{equation*}
\tilde{f} (g) := \sum\limits_{v\in \Z^n} {} f(g v).
\end{equation*}

Lemma 3.1 in \cite{S2} implies that
\begin{equation}
\label{eq:tld f const alpha} \tilde{f}(g) < c \alpha(g),
\end{equation}
where $c=c(f)$ is an explicit constant constant
\begin{equation*}
c(f) = c_0 M \max(1, R^n),
\end{equation*}
for some constant $c_0 = c_0 (n)$, independent on f. In section
\ref{subsec:mem const exm} we prove a stronger result, assuming
some additional information about the support of $f$.

Let $Q_0$ be a quadratic form defined by
\begin{equation*}
Q_0(\vec{v}) = 2v_1 v_n + \sum\limits_{i=2}^{p} v_i ^2 -
\sum\limits_{i=p+1}^{n-1} v_i ^2.
\end{equation*}
Since
\begin{equation*}
v_1 v_n = \frac{(v_1+v_n)^2 - (v_1-v_n)^2}{2},
\end{equation*}
$Q_0$ is of signature $p,\, q$. Obviously, $G:=SL(n,\R)$ acts on
the space of quadratic forms of signature $(p,\,q)$, and
discriminant $\pm 1$, $\OO=\OO(p,\,q)$ by:
\begin{equation*}
Q^{g} (v) := Q(gv).
\end{equation*}
Moreover, by the well known classification of quadratic forms,
$\OO$ is the orbit of $Q_0$ under this action.

In our case the signature is $(p,\, q) = (2,\, 2)$ and $n = 4$. We
fix an element $h_1 \in G$ with $Q^{h_1} = Q_1$, where $Q_1$ is
given by \eqref{eq:quad form def}. There exists a constant $\tau >
0$, such that for every $v\in \R^4$,
\begin{equation}
\label{eq:tau prop} \tau ^{-1} \| v \| \le \| h_1 v \| \le \tau \|
v \|.
\end{equation}
We may assume, with no loss of generality that $\tau \ge 1$.

Let $H:=Stab_{Q_0} (G)$. Then the natural mophism $H \backslash G
\rightarrow \OO(p,q)$ is a homeomorphism. Define a $1$-parameter
family $a_t\in G$ by:
\begin{equation*}
a_t e_i = \begin{cases} e^{-t} e_1, \quad &i=1 \\ e_i,\quad
&i=2,\ldots,\, n-1
\\ e^t e_n,\, &i=n
\end{cases}.
\end{equation*}
Clearly, $a_t\in H$. Furthermore, let $\hat{K}$ be the subgroup of
$G$ consisting of orthogonal matrices, and denote $K :=
H\cap\hat{K}$.

Let $(a,\, b) \in \R^2$ be given and let $Q:\R^n\rightarrow\R$ be
any quadratic form. The object of our interest is:
\begin{equation*}
V_{(a,\, b)} (\Z) = V_{(a,\, b)} ^ Q (\Z) = \{x\in\Z^n :\: a <
Q(x) < b\}.
\end{equation*}

Theorem 2.3 states, in our case:
\begin{theorem}[Theorem 2.3 from ~\cite{EMM}]
\label{thm:mem thm 2.3} Let $\Omega = \{ v\in\R^4 |\: \| v \| <
\nu(v/\| v \|) \}$, where $\nu$ is a nonnegative continuous
function on $S^{3}$. Then we have:
\begin{equation*}
\# V_{(a,\,b)} ^ {Q_1} (\Z) \cap T\Omega < c T^{2} \log{T},
\end{equation*}
where the constant $c$ depends only on $(a,\, b)$.
\end{theorem}

The proof of theorem \ref{thm:mem thm 2.3} relies on theorem 3.3
from ~\cite{EMM}, and we give here a particular case of this
theorem
\begin{theorem}[Theorem 3.3 from ~\cite{EMM}]
\label{thm:mem thm 3.3} For any (fixed) lattice $\Delta$ in
$\R^4$,
\begin{equation*}
\sup\limits_{t>1} \frac{1}{t} \int\limits_{K} \alpha(a_t k \Delta)
dm(k) < \infty,
\end{equation*}
where the upper bound is universal.
\end{theorem}

\subsection{Outline of the proof of theorem \ref{thm:mem thm 2.3}:}

\paragraph{Step 1:}
Define
\begin{equation}
\label{eq:Jf def} J_f (r, \zeta) = \frac{1}{r^{2}} \int
\limits_{\R^{2}} f(r, x_2,\, x_{3}, \, x_4) dx_2 dx_{3},
\end{equation}
where
\begin{equation*}
x_4 = \frac{\zeta - x_2 ^2 + x_3 ^2}{2r}
\end{equation*}

Lemma 3.6 in ~\cite{EMM} states that $J_f$ is approximable by
means of an integral over the compact subgroup K. More precisely,
there is some constant $C>0$, such that for every $\epsilon > 0$,
\begin{equation}
\label{eq:Jf int approx} \bigg| C \cdot e^{2t} \int\limits_{K}{}
f(a_t kv) \nu (k^{-1} e_1) dm(k) - J_f\big(\| v \| e^{-t},\, Q_0
(v)\big) \nu (\frac{v}{\| v \|}) \bigg| < \epsilon
\end{equation}
with $e^t,\, \| v \| > T_0$ for some $T_0 > 0$.

\paragraph{Step 2:}
Choose a continuous nonnegative function $f$ on $\R_{+} ^4 = \{
x_1 > 0\}$ which vanishes outside a compact set so that
\begin{equation*}
J_f (r,\zeta) \ge 1+\epsilon
\end{equation*}
on $[\tau^{-1},\, 2\tau]\times [a,\, b]$. We will show later, how
one can choose $f$.

\paragraph{Step 3:}
Denote $T=e^t$, and suppose that $T\le \| v \| \le 2T$ and $a\le
Q_0 (h_1 v) \le b$. Then by \eqref{eq:tau prop}, $J_f\big(\| h_1 v
\| T^{-1}, Q_0 (h_1 v)\big) \ge 1+\epsilon$, and by \eqref{eq:Jf
int approx}, for a sufficiently large $t$,
\begin{equation}
\label{eq:V and int} C \cdot T^{2} \int\limits_{K} f(a_t k h_1 v)
dm(k) \ge 1,
\end{equation}
for $T\le \| v \| \le 2T$ and
\begin{equation}
\label{eq:Qx in int} a\le Q_0^x (v) \le b.
\end{equation}

\paragraph{Step 4:}
Summing \eqref{eq:V and int} over all $v\in \Z^4$ with
\eqref{eq:Qx in int} and $T\le \| v \| \le 2T$, we obtain:
\begin{equation}
\label{eq:V bnd tld}
\begin{split}
\# V_{(a,\,b)} (\Z) \cap [T,\,2T]S^{3} &\le
\sum\limits_{v\in \Z^n} {} C \cdot T^{2} \int\limits_{K} f(a_t k h_1 v) dm(k) \\
&=  C \cdot T^{2} \int\limits_{K} \tilde{f} (a_t k h_1) dm(k)
\end{split}
\end{equation}
using the nonnegativity of $f$.

\paragraph{Step 5:}
By \eqref{eq:tld f const alpha}, \eqref{eq:V bnd tld} is
\begin{equation*}
\le C\cdot c(f) \cdot T^{2} \int\limits_{K} \alpha (a_t k h_1)
dm(k).
\end{equation*}

\paragraph{Step 6:}
The result of theorem 2.3 is obtained by using theorem
\ref{thm:mem thm 3.3} on the last expression.

\subsection{$\delta$-dependency:}
\label{subsec:mem const exm}

In this section we assume that $(a,\, b) = (0,\, \delta)$, which
suits the definition of the set $A(R,\, \delta)$,
\eqref{eq:A(R,delta) def}. One should notice that there only $3$
$\delta$-dependent steps:

$\bullet$ Choosing $f$ in step 2, such that $J_f \ge 1+\epsilon$
on $[\tau^{-1},\, 2\tau]\times [0,\, \delta]$. We will construct a
family of functions $f_{\delta}$ with an universal bound
$|f_{\delta}| \le M$, such that $f_{\delta}$ vanishes outside of a
compact set which is only slightly larger than
\begin{equation}
\label{eq:V delta def} V(\delta) = [\tau^{-1},\, 2\tau] \times
[-1,\, -1]^{2} \times [0,\, \frac{\delta \tau}{2}].
\end{equation}
This is done in section \ref{subsubsec:choose f}.

$\bullet$ The dependency of $T_0$ of step 3, so that the usage of
lemma 3.6 in ~\cite{EMM} is legitimate. For this purpose we will
have to examine the proof of this lemma. This is done in section
\ref{subsubsec:T_0 dep}.

$\bullet$ The constant $c$ in \eqref{eq:tld f const alpha}. We
would like to establish a {\em linear} dependency on $\delta$.
This is straightforward, once we are able to control the number of
integral points in a domain defined by \eqref{eq:V delta def}.
This is done in section \ref{subsubsec:int V pnt cntrl}.

\subsubsection{Choosing $f_{\delta}$:}
\label{subsubsec:choose f}

\paragraph{Notation:} For a set $U\subset \R^n$, and $\epsilon > 0$, denote
\begin{equation*}
U_{\epsilon} := \{ x\in \R^n:\: \max\limits_{1\le i \le n}
|x_i-y_i| \le \epsilon, \,\text{for some } y\in U\}.
\end{equation*}

Choose a nonnegative continuous function $f_0$, on $\R^4_{+}$,
which vanishes outside a compact set, such that its support,
$E_{f_0}$, slightly exceeds the set $V(1)$. More precisely, $V(1)
\subset E_{f_0} \subset V(1)_{\delta_0}$ for some $\delta_0 > 0$.
By the uniform continuity of $f$, there are $\epsilon_0,\,
\delta_0 > 0$, such that if $\max\limits_{1\le i \le 4} |x_i - x_i
^0| \le \delta_0$, then $f(x) > \epsilon_0$, for every  $x^0 =
(x_1^0,0,\,0,\, x_4^0)\in V(1)$.

Thus for $(r,\,\zeta) \in [\tau^{-1},\, 2\tau]\times [0,\,
\delta]$, the contribution of $[-\delta_0,\, \delta_0] ^ {2}$ to
$J_{f_0}$ is $\ge \epsilon_0 \cdot (2\delta_0)^{2}$. Multiplying
$f_0$ by a suitable factor, and by the linearity of $J_{f_0}$, we
may assume that this contribution is at least $1+\epsilon$.

Now define $f_{\delta} (x_1,\,\ldots,\, x_4) := f_0 (x_1,\,x_2,\,
x_{3},\, \frac{x_4}{\delta})$. We have for $\delta \ge 1$
\begin{equation*}
\frac{\zeta - x_2^2 + x_3^2}{2r\delta} = \frac{\zeta / 2r}{\delta}
- \frac{(x_2 /\sqrt{\delta})^2}{2r} + \frac{(x_3
/\sqrt{\delta})^2}{2r}.
\end{equation*}
Thus for $\delta \ge 1$, if $(r,\, \zeta) \in [\tau^{-1},\,
2\tau]\times [0,\, \delta]$ and for $i=2,\, 3$, $|x_i| < \delta_0
$, $f_{\delta}$ satisfies:
\begin{equation*}
f_{\delta} (r,\, x_2,\, x_{3},\, x_4) > \epsilon_0,
\end{equation*}
and therefore the contribution of this domain to $J_{f_{\delta}}$
is
\begin{equation*}
\ge \epsilon_0 (2\delta)^2 \ge 1+\epsilon
\end{equation*}
by our assumption.

By the construction, the family $\{ f_{\delta} \} $ has a
universal upper bound $M$ which is the one of $f_0$.

\subsubsection{How large is $T_0$}
\label{subsubsec:T_0 dep}

The proof of lemma 3.6 from ~\cite{EMM} works well along the same
lines, as long as
\begin{equation}
\label{eq:x in sup} f(a_t x) \ne 0
\end{equation}
implies that for $t\rightarrow \infty$, $x / \| x \|$ converges to
$e_1 = (1,\, 0,\, 0,\, 0)$. Now, since $a_t$ preserves $x_1 x_4$,
\eqref{eq:x in sup} implies for the particular choice of
$f=f_{\delta}$ in section \ref{subsubsec:choose f}:
\begin{equation*}
|x_1 x_4 | = O(\delta) ; \quad x_1 \gg T.
\end{equation*}
Thus
\begin{equation*}
\| x \| = x_1 + O(\frac{\delta}{T}) + O(1),
\end{equation*}
and so, as long as $\delta = o(T)$, $x / \| x \|$ indeed converges
to $e_1$.

\subsubsection{Bounding integral points in $V_{\delta}$:}
\label{subsubsec:int V pnt cntrl}

\begin{lemma}
\label{lem:V delta lat bnd} Let $V(\delta)$ defined by
\begin{equation}
V(\delta) = [\tau^{-1},\, 2\tau] \times [-1,\, -1]^{n-2} \times
[0,\, \frac{\delta \beta}{2}].
\end{equation}
for some constant $\tau$ and $n \ge 3$. Let $g\in SL(n,\,\R)$ and
denote
\begin{equation*}
N(g,\, \delta) := \# V(\delta) \cap g\Z^n.
\end{equation*}
Then for $\delta \ge 1$,
\begin{equation*}
\bigg| N(g,\, \delta) - \frac{2^{n-2} (2\tau-\tau^{-1})\delta
}{\det {g}} \bigg| \le c_5 \delta \sum\limits_{i=1}^{n-1}
\frac{1}{vol(L_i / (g\Z^n \cap L_i)}
\end{equation*}
for some $g$-rational subspaces $L_i$ of $\R^4$ of dimension $i$,
where $c_5 = c_5(n)$ depends only on $n$.
\end{lemma}

A direct consequence of lemma \ref{lem:V delta lat bnd} is the
following
\begin{corollary}
Let $f:\R^n\rightarrow \R$ be a nonnegative function which
vanishes outside a compact set $E$. Suppose that $E \subset
V_{\epsilon} (\delta)$ for some $\epsilon > 0$. Then for $\delta
\ge 1$, \eqref{eq:tld f const alpha} is satisfied with
\begin{equation*}
c(f) = c_3 \cdot M \delta,
\end{equation*}
where the constant $c_3$ depends on $n$ only.
\end{corollary}

In order to prove lemma \ref{lem:V delta lat bnd}, we shall need
the following:

\begin{lemma}
Let $\Lambda\subset\R^n$ be a $m$-dimensional lattice, and let
\begin{equation}
\label{eq:A delta def} A_{t} = \begin{pmatrix} 1 \\ & 1 \\ & &
\ddots \\ & & & t \end{pmatrix}
\end{equation}
an $n$-dimensional linear transformation. Then for $t > 0$ we have
\begin{equation}
\label{eq:sub lin det} \det{A_{t}\Lambda} \le t \det{\Lambda}.
\end{equation}
\end{lemma}
\begin{proof}
We may assume that $m < n$, since if $m=n$, we obviously have an
equality. Let $v_1,\,\ldots,\, v_m$ the basis of $\Lambda$ and
denote for every $i$, $u_i\in\R^{n-1}$ the vector, which consists
of first $n-1$ coordinates of $v_i$. Also, let $x_i\in\R$ be the
last coordinate of $v_i$. By switching vectors, if necessary, we
may assume $x_1 \ne 0$. We consider the function
\begin{equation*}
f(t) := (\det{A_{t}\Lambda}) ^2,
\end{equation*}
as a function of $t\in\R$.

Obviously,
\begin{equation*}
f(t) = \det \big( <u_i,\, u_j> + x_i x_j t^2 \big)_{1 \le i,\, j
\le m}.
\end{equation*}
Substracting $\frac{x_i}{x_1}$ times the first row from any other,
we obtain:
\begin{equation*}
f(t) = \begin{vmatrix}  <u_1,\, u_j> +x_1 x_j t^2\\ <u_2,\, u_j >
-
\frac{x_2}{x_1} <u_1,\, u_j> \\  \vdots \\
<u_m,\, u_j > - \frac{x_m}{x_1} <u_1,\, u_j> \end{vmatrix},
\end{equation*}
and by the multilinearity property of the determinant, $f$ is a
linear function of $t^2$. Write
\begin{equation*}
f(t) = a(t^2-1)+b t^2.
\end{equation*}
Thus
\begin{equation*}
b = f(1) ;  \quad a = -f(0),
\end{equation*}
and so $b = \det{\Lambda}$, and $a = -\det < u_i,\, u_j> \le 0$,
being minus the determinant of a Gram matrix. Therefore,
\begin{equation*}
(\det{A_{t}\Lambda})^2 - t^2 \det{\Lambda} = a (t^2-1) \le 0
\end{equation*}
for $t \ge 1$, implying \eqref{eq:A delta def}.
\end{proof}

\begin{proof}[Proof of lemma \ref{lem:V delta lat bnd}]
We will prove the lemma, assuming $\beta=2$. However, it implies
the result of the lemma for any $\beta$, affecting only $c_5$. Let
$\delta > 0$. Trivially,
\begin{equation*}
N(g,\, \delta) = N(g_0,\, 1),
\end{equation*}
where $g_0 = A_{\delta} ^ {-1} g$ with $A_{\delta}$ given by
\eqref{eq:A delta def}. Let
$\lambda_1\le\lambda_2\le\ldots\le\lambda_n$ be the successive
minima of $g_0$, and pick linearly independent lattice points
$v_1,\,\ldots,\, v_n$ with $\| v_i \| = \lambda_i$. Denote $M_i$
the linear space spanned by $v_1,\,\ldots,\, v_i$ and the lattice
$\Lambda_i = g_0\Z^n \cap M_i$.

First, assume that $\lambda_n \le \sqrt{\tau^2 + (n-1)} =: r$.
Now, by Gauss' argument,
\begin{equation*}
\bigg| N(g_0,\, 1) - \frac{2^{n-1} (2\tau-\tau^{-1}) \delta}{\det
{g}} \bigg| \le \frac{1}{\det{g_0}} vol (\Sigma),
\end{equation*}
where
\begin{equation*}
\Sigma := \{ x:\: dist(x,\, \partial V(1))\le n\lambda_n \}.
\end{equation*}
Now, for $\lambda_n \le r$,
\begin{equation*}
vol(\Sigma) \ll \lambda_n,
\end{equation*}
where the constant implied in the $``\ll``$-notation depends on
$n$ only (this is obvious for $\lambda_n \le \frac{1}{2n}$, and
trivial otherwise, since for $\lambda_n \le r$, $vol(\Sigma) =
O(1)$). Thus,
\begin{equation*}
\begin{split}
\bigg| &N(g_0,\, 1) - \frac{2^{n-1} (2\tau-\tau^{-1}) \delta}{\det
{g}} \bigg| \ll \frac{\lambda_n}{\det{g_0}} \ll
\frac{1}{\det{\Lambda_{n-1}}} \\&= \frac{1}{vol(M_{n-1} / M_{n-1}
\cap g_0 \Z^n)} \le \frac{\delta}{vol(A_{\delta} M_{n-1} /
A_{\delta} M_{n-1} \cap g \Z^n )}
\end{split}
\end{equation*}

Next, suppose that $\lambda_n > r$. Then,
\begin{equation*}
V(\delta) \cap g_0 \Z^n \subset V(\delta) \cap \Lambda_{n-1}.
\end{equation*}
Thus, by the induction hypothesis, the number of such points is:
\begin{equation*}
\begin{split}
\le &c_4 \sum\limits_{i=0}^{k-1} \frac{1}{\det (\Lambda_i)} =
\sum\limits_{i=0}^{k-1} \frac{1}{vol(M_i / M_i \cap g_0\Z^n)}  \\
&\le  \delta \sum\limits_{i=0}^{k-1} \frac{1}{vol(A_{\delta} M_i /
A_{\delta} \M_i \cap g\Z^n)}.
\end{split}
\end{equation*}
Since $\lambda_n > r$, we have
\begin{equation*}
\frac{1}{\det{g}} = \frac{1}{\lambda_n} \frac{1}{\det{g} /
\lambda_n} \ll \frac{1}{\det{g} / \lambda_n} \ll
\frac{1}{\lambda_1\cdot\ldots\cdot\lambda_{n-1}},
\end{equation*}
and we're done by defining $L_i := A_{\delta} M_i$.
\end{proof}

\section{Unsmoothing}
\label{sec:unsmth}

\subsection{An asymptotic formula for $N_{\Lambda}$}
\label{sec:asym for N_Lambda}

We need  an asymptotic formula for the {\em sharp} counting
function $N_{\Lambda}$. Unlike the case of the standard lattice,
$\Z^2$, in order to have a good control over the error terms we
should use some Diophantine properties of the lattice we are
working with. We adapt the following notations:

Let $\Lambda$ be a lattice and $t > 0$ a real variable. Denote the
set of squared norms of $\Lambda$ by
\begin{equation*}
SN_{\Lambda} = \{ |\vec{n}|^2:\: n\in\Lambda \}.
\end{equation*}
Suppose we have a function $\delta_\Lambda
:SN_{\Lambda}\rightarrow \R$, such that given $\vec{k} \in \Lambda
$, there are \underline{no} vectors $\vec{n} \in \Lambda$ with $0
< | |\vec{n}|^2 - |\vec{k}|^2 | < \delta_{\Lambda} ( |\vec{k} | ^2
)$. That is,
\begin{equation*}
\Lambda \cap \{ \vec{n}\in\Lambda :\: |\vec{k}|^2 -
\delta_{\Lambda} ( |\vec{k} | ^2 ) < | \vec{n} |^2 < |\vec{k}|^2 +
\delta_{\Lambda} ( |\vec{k} | ^2 ) \} = A_{|\vec{k}|},
\end{equation*}
where
\begin{equation*}
A_{y} := \{ \vec{n}\in\Lambda :\: |\vec{n}| = y\}.
\end{equation*}
Extend $\delta_{\Lambda}$ to $\R$ by defining $\delta_{\Lambda}
(x) := \delta_{\Lambda} (|\vec{k}|^2 )$, where $\vec{k} \in
\Lambda$ minimizes $|x-|\vec{k} |^2|$ (in the case there is any
ambiguity, that is if $x = \frac{|\vec{n_1}|^2+|\vec{n_2}|^2}{2}$
for vectors $\vec{n_1},\,\vec{n_2} \in \Lambda$ with consecutive
increasing norms, choose $\vec{k}:= \vec{n_1}$). We have the
following lemma:

\begin{lemma}
\label{lem:asym cnt shrp} For every $a>0,\,c>1$,
\begin{equation*}
\begin{split}
N_{\Lambda } (t) &= \frac{\pi }{\beta} t^2 - \frac{\sqrt{t}}{\beta
\pi } \sum\limits_{\substack{\vec{k} \in \Lambda ^ {*} \setminus
\{ 0 \}
\\ |\vec{k} | \le \sqrt{N} }} \frac{\cos \big( 2\pi t |\vec{k} | +
\frac{\pi}{4} \big) } {|\vec{k}| ^ \frac{3}{2}} +
O(N^{a}) \\
&+ O \bigg(\frac{t^{2c-1}}{\sqrt{N}} \bigg) + O
\bigg(\frac{t}{\sqrt{N}}
\cdot \big( \log t + \log (\delta_{\Lambda} (t^2) \big) \bigg) \\
&+ O \bigg( \log{N} + \log (\delta_{\Lambda ^ *} (t^2) ) \bigg)
\end{split}
\end{equation*}
\end{lemma}
As a typical example of such a function, $\delta_{\Lambda }$, for
$\Lambda = \langle 1, \, \alpha + i \beta \rangle$, with a
Diophantine $(\alpha,\, \alpha^2,\, \gamma^2)$, we may choose
$\delta_{\Lambda } (y) = \frac{c}{y^{K}}$, where $c$ is a
constant. In this example, if $\Lambda\ni\vec{k} = (a,b)$, then by
lemma \ref{lem:no mult}, $A_{|\vec{k}|} = \pm(a,\, b )$, provided
that $\gamma$ is irrational.

The proof of this lemma is essentially the same as the one of
lemma 5.1 in ~\cite{W}, starting from
\begin{equation*}
\ZZ_{\Lambda} (s) := \frac{1}{4} \sum\limits_{\vec{k} \in \Lambda
\setminus {0}} \frac{1}{|\vec{k}|^{2s}}
 = \sum\limits_{(m,\, n ) \in \Z_{+}^{2} \setminus {0}} \frac{1}{\big(
(m+n\alpha)^2+ (\beta n)^2 \big)^{s}}
\end{equation*}

\begin{proposition}
\label{prop:bnd 2nd mom diff} Let a lattice $\Lambda = \langle 1,
\, \alpha+i\beta \rangle$ with a Diophantine triple of numbers
$(\alpha^2, \alpha\beta, \, \beta^2)$ be given. Suppose that $L
\rightarrow \infty$ as $T\rightarrow \infty$ and choose $M$, such
that $L/\sqrt{M} \rightarrow 0$, but $M = O \big( T^\delta \big)$
for every $\delta > 0$ as $T \rightarrow \infty$. Suppose
furthermore, that $M = O (L^{s_0})$ for some (fixed) $s_0 > 0$.
Then
\begin{equation*}
\Bigg\langle \bigg| S_{\Lambda} (t,\, \rho ) -
\tilde{S}_{\Lambda, \, M, \, L} (t) \bigg| ^ {2} \Bigg\rangle \ll
\frac{1}{\sqrt{M}}
\end{equation*}
\end{proposition}

The proof of proposition \ref{prop:bnd 2nd mom diff} proceeds
along the same lines as the one of proposition 6.1 in ~\cite{W},
using again an asymptotic formula for the sharp counting function,
given by lemma \ref{lem:asym cnt shrp}. The only difference is
that here we use proposition \ref{prop:S_n bnd blh} rather than
lemma 6.2 from ~\cite{W}.

Once we have proposition \ref{prop:bnd 2nd mom diff} in our hands,
the proof of our main result, namely, theorem \ref{thm:norm dist}
proceeds along the same lines as the one of theorem 1.1 in
~\cite{W}.

\paragraph{Acknowledgement.}
This work was supported in part by the EC TMR network
\textit{Mathematical Aspects of Quantum Chaos}, EC-contract no
HPRN-CT-2000-00103 and the Israel Science Foundation founded by
the Israel Academy of Sciences and Humanities. This work was
carried out as part of the author's PHD thesis at Tel Aviv
University, under the supervision of prof. Ze\'{e}v Rudnick. The
author wishes to thank Alex Eskin for his help. A substantial part
of this work was done during the author's visit to the university
of Bristol.

\end{document}